\documentclass[12pt]{article}
\usepackage{amscd}
\usepackage{amsmath,amsfonts,amssymb,amscd}
\usepackage{indentfirst,graphics,epsfig}
\input{epsf}

\usepackage{enumerate}

\usepackage{enumerate}

\usepackage{lineno}

\makeatletter \@addtoreset{figure} {section}\makeatother
\makeatletter
\long\def\@makecaption#1#2{%
   \vskip 10\p@
   \setbox\@tempboxa\hbox{{#1}\ \ #2}%
   \ifdim \wd\@tempboxa >\hsize
       {#1}\ \ #2\par
   \else
       \hbox to\hsize{\hfil\box\@tempboxa\hfil}%
   \fi}
\makeatother

\newtheorem{thm}{Theorem}
\newtheorem{defi}[thm]{Definition}
\newtheorem{cor}[thm]{Corollary}
\newtheorem{lem}[thm]{Lemma}

\newtheorem{con}[thm]{Conjecture}

\newcommand{\qed}{{\hfill\rule{3pt}{7pt}}}
\def\pf{\noindent {\it Proof.} }

\def\qed{\hfill \rule{4pt}{7pt}}
\def\pf{\noindent {\it Proof.} }

\title{\bf Note for Nikiforov's two conjectures\\
on the energy of trees\footnote{Supported by NSFC No.10831001,
PCSIRT and the ``973" program.}}

\author{ \small Xueliang Li,~~Jianxi Liu\\
\small Center for Combinatorics and LPMC-TJKLC\\
\small Nankai University, Tianjin 300071, P.R. China\\
\small Email: lxl@nankai.edu.cn}
\begin{document}
\date{}
\maketitle
\begin{abstract}
{\small The energy $E$ of a graph is defined to be the sum of the
absolute values of its eigenvalues. Nikiforov in {\it ``V.
Nikiforov, The energy of $C_4$-free graphs of bounded degree, Lin.
Algebra Appl. 428(2008), 2569--2573"} proposed two conjectures
concerning the energy of trees with maximum degree $\Delta\leq 3$.
In this short note, we show that both conjectures are true. }\\[3mm]
{\bf Key words:} energy of a graph, conjecture, tree\\[2mm]
{\bf AMS Subject Classification:} 05C50, 05C90, 15A18, 92E10
\end{abstract}

\vskip1.5cm

Let $G$ be a graph on $n$ vertices and $\lambda_1,\ \lambda_2,
\cdots, \lambda_n$ be the eigenvalues of its adjacency matrix. The
value $E(G) = |\lambda_1| + \cdots + |\lambda_n|$ is defined as the
energy of $G$, which has been studied intensively, see \cite{G1,
GLZ} for a survey.

In \cite{N}, Nikiforov proposed two conjectures on the energy of
trees. In order to state and prove them, we need the following
notations and terminology.

\begin{figure}[!hbpt]
\begin{center}
\includegraphics[scale=0.70000]{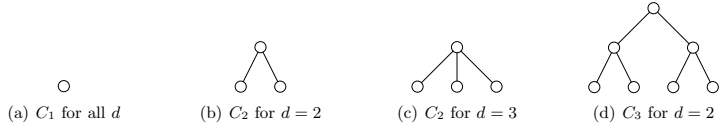}
\end{center}\vspace{-15pt}
\caption{Some small complete $d$-ary trees.}\label{figure1}
\end{figure}

The {\it complete $d$-ary tree} of height $h-1$ is denoted by
$C_h$, which is built up inductively as follows: $C_1$ is a single
vertex and $C_h$ has $d$ branches $C_{h-1},\cdots, C_{h-1}$. See
Figure \ref{figure1} for examples.
It is convenient to set $C_0$ as the empty graph.

Let $\mathcal{T}_{n,d}$ be the set of all trees with $n$ vertices
and maximum degree $d + 1$. We define a special tree $T_{n,d}^*$ as
follows (see also \cite{HW}):
\begin{defi} $T_{n,d}^*$ is the tree with $n$ vertices that can be decomposed
as in Figure \ref{figure2} \vspace{-8pt}
\begin{figure*}[!hbpt]
\begin{center}
\includegraphics[scale=0.70000]{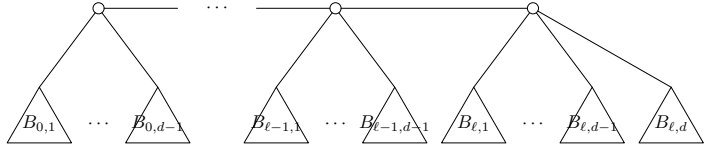}
\end{center}\vspace{-15pt}
\caption{Tree $T_{n,d}^*$.}\label{figure2}
\end{figure*}

\noindent with $B_{k,1}, \cdots, B_{k,d-1}\in \{C_k,C_{k+2}\}$ for
$0 \leq k < l$ and either $B_{l,1} =\cdots = B_{l,d} = C_{l-1}$ or
$B_{l,1} =\cdots = B_{l,d} = C_l$ or $B_{l,1},\ \cdots,\ B_{l,d}
\in \{C_{l},C_{l+1},C_{l+2}\}$, where at least two of $B_{l,1},
\cdots, B_{l,d}$ equal $C_{l+1}$. This representation is unique,
and one has the ``digital expansion''
\begin{align} (d - 1)n + 1 =
\sum\limits_{k=0}^l a_k d^k, \label{3}
\end{align}
where $a_k =
(d-1)(1+(d+1)r_k)$ and $0 \leq r_k \leq d-1$ is the number of
$B_{k,i}$ that are isomorphic to $C_{k+2}$ for $k < l$, and
\begin{align*}
&\bullet a_l = 1 \ \ \mbox{if}\ \ B_{l,1} = \cdots = B_{l,d} =
C_{l-1},\\
&\bullet a_l = d \ \ \mbox{if}\ \ B_{l,1}= \cdots = B_{l,d} = C_l,\\
&\bullet \mbox{or otherwise}\ a_l = d + (d - 1)q_l + (d^2 - 1)r_l,\
\mbox{where}\ q_l\geq 2\ \mbox{is the number of}\ B_{l,i}\ \\
&\mbox{that are isomorphic to}\ C_{l+1}\  \mbox{and}\ r_l \
\mbox{is the number of}\ B_{l,i}\ \mbox{that are isomorphic to}\
C_{l+2}.
\end{align*}
\end{defi}

Let $\mathcal{B}_n$ denote the tree constructed by taking three
disjoint copies of the complete $2$-ary tree of height $h-1$,
i.e., $C_n$, and joining an additional vertex to their roots
(i.e., vertices of height zero). In the end of \cite{N}, Nikiforov
formulated two conjectures as follows:

\begin{con}\label{con1} The limit\[ c=\lim\limits_{n\rightarrow \infty} \frac{E(\mathcal{B}_n)}{3\cdot
2^{n+1}-2}\] exists and $c>1$.
\end{con}

\begin{con}\label{con2}
Let $\epsilon > 0$. If $T$ is a sufficiently large tree with
$\Delta(T) \leq 3$, then $E(T) \geq (c -\epsilon)|T|$.
\end{con}

Nikiforov mentioned that empirical data given in \cite{G2} seem to
corroborate these conjectures, but apparently new techniques are
necessary to prove or disprove them. We will give confirmative
proofs for both Conjecture \ref{con1} and Conjecture \ref{con2}.

We first state two known lemmas from \cite{HW}, which will be needed
in the sequel.

\begin{lem}\cite{HW}.\label{lem1} Let $n$ and $d$ be positive integers. Then $T_{n,d}^*$ is the
unique (up to isomorphism) tree in $\mathcal{T}_{n,d}$ that
minimizes the energy.
\end{lem}

\begin{lem}\cite{HW}.\label{lem2} The energy of $T_{n,d}^*$ is asymptotically
\[ E(T_{n,d}^*) =\alpha_d \cdot n+O(\ln n),\]
where
\begin{align}\label{4}
\alpha_d=2\sqrt{d} (d-1)^2 \Bigg{(} \sum\limits_{j\geq 1\atop
j\equiv 0\ (mod\ 2)} d^{-j} \left( \cot \frac{\pi}{2j}-1 \right)+
\sum\limits_{j\geq 1\atop j\equiv 1\ (mod\ 2)} d^{-j} \left(\csc
\frac{\pi}{2j}-1 \right) \Bigg{)}
\end{align}
is a constant that only depends on $d$.
\end{lem}
\begin{table}[!hbpt]
\begin{center}
\includegraphics[scale=1.0000]{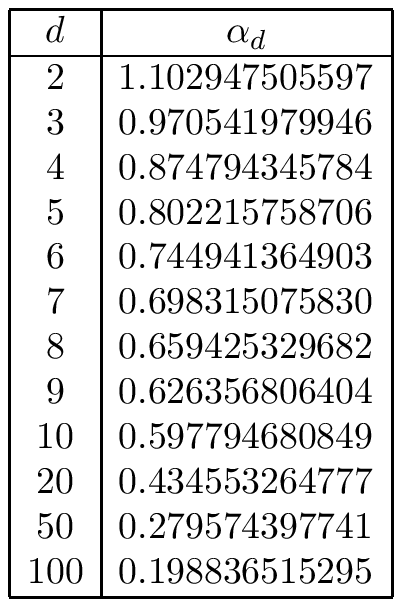}
\caption{Some numerical values for the constant
$\alpha_d$.}\label{table1}
\end{center}\vspace{-15pt}
\end{table}

With the above two lemmas, the two conjectures can be proved very
easily as follows.

\begin{thm} The limit\[ c=\lim\limits_{n\rightarrow \infty} \frac{E(\mathcal{B}_n)}{3\cdot
2^{n+1}-2}\] exists and $c>1$.
\end{thm}
\pf We just need to notice that $\mathcal{B}_n$ is exactly the tree
$T_{3\cdot 2^{n+1}-2,2}^*$ with $l=n,\ B_{k,1}=C_k$ for $0 \leq k <
l,\ B_{l,1} = B_{l,2} = C_n$. Therefore, by Lemma \ref{lem2} and
Table \ref{table1} we have $$ \lim\limits_{n\rightarrow \infty}
\Bigg{(}\frac{E(\mathcal{B}_n)}{3\cdot 2^{n+1}-2}
\Bigg{)}=\lim\limits_{n\rightarrow \infty} \Bigg{(} \alpha_2+
\frac{O(\ln (3\cdot 2^{n+1}-2))}{{3\cdot 2^{n+1}-2}}
\Bigg{)}=\alpha_2>1.$$\qed

In fact, from Lemmas \ref{lem1} and \ref{lem2} we have that for any
$T\in \mathcal{T}_{n,d}$, $$E(T)\geq E(T_{n,d}^*)= \alpha_d\cdot n +
O(\ln n).$$ Therefore, we obtain
\begin{thm}\label{thm2}
Let $\epsilon > 0$. If $T$ is a sufficiently large tree with
$\Delta(T) = d+1$, then $E(T) \geq (\alpha_d-\epsilon)|T|$, where
$\alpha_d$ is given in Equ.\eqref{4}.
\end{thm}

Letting $d=2$, we get \begin{cor} Let $\epsilon > 0$. If $T$ is a
sufficiently large tree with $\Delta(T) = 3$, then $E(T) \geq
(\alpha_2-\epsilon)|T|$, where $\alpha_2$ is given in Equ.\eqref{4}.
\end{cor}

Recall that a hypoenergetic graph of order $n$ is such that $E(G)<
n$, whereas it is strongly hypoenergetic if $E(G)< n-1$. We have the
following easy remarks:

\noindent {\bf Remark 1:} From Lemma \ref{lem2} and Table
\ref{table1}, one can see that there is neither strongly
hypoenergetic tree nor hypoenergetic tree of order $n$ and maximum
degree $\Delta$ for $\Delta \leq 3$ and any suitable large $n$.

\noindent {\bf Remark 2:} From Lemma \ref{lem2} and Table
\ref{table1}, one can also see that there are both hypoenergetic
trees and strongly hypoenergetic trees of order $n$ and maximum
degree $\Delta$ for $\Delta \geq 4$ and any suitable large $n$.

\end{document}